\documentclass[12pt]{article}
\nonstopmode
\usepackage{graphicx,xspace,colortbl}
\usepackage{amsmath,amsthm,amsfonts}
\usepackage{fancybox,float}
\restylefloat{figure}
\RequirePackage[colorlinks,citecolor=blue,urlcolor=blue,linkcolor=blue]{hyperref}
\hypersetup{
colorlinks = true,
citecolor=blue,
urlcolor=blue,
linkcolor=blue,
pdfauthor = {A. Kuznetsov, J.C. Pardo},
pdfkeywords = {Hypergeometric Levy processes, Lamperti-stable processes, exponential functional, double gamma function, Lamperti transformation, extrema of stable processes},
pdftitle = {Fluctuations of stable processes and exponential functionals of hypergeometric  Levy processes},
pdfpagemode = UseNone
}
    \def\qed{\hfill$\sqcap\kern-8.0pt\hbox{$\sqcup$}$\\}
    \def\beq{\begin{eqnarray}}
    \def\eeq{\end{eqnarray}}
    \def\beqq{\begin{eqnarray*}}
    \def\eeqq{\end{eqnarray*}}

    \def\re{\textnormal {Re}}
    \def\im{\textnormal {Im}}
    \def\p{{\mathbb P}}
    \def\e{{\mathbb E}}
    \def\r{{\mathbb R}}
    \def\mm{{\mathcal M}}
    
    \def\c{{\mathbb C}}

    \def\d{{\textnormal d}}
    \def\i{{\textnormal i}}

  \def\txt{\textstyle}
	\newtheorem{theorem}{Theorem}
	\newtheorem{lemma}{Lemma}
	\newtheorem{proposition}{Proposition}
	\newtheorem{corollary}{Corollary}
	\newtheorem{definition}{Definition}

\title{
\textbf{Fluctuations of stable processes and exponential functionals of hypergeometric  L\'evy processes}
}
\author{
\textbf{
A. Kuznetsov
\footnote{Department of Mathematics and Statistics, 
York University, 
4700 Keele Street, 
Toronto, Ontario, 
M3J 1P3, Canada. Email: kuznetsov@mathstat.yorku.ca.
Research supported by the
Natural Sciences and Engineering Research Council of Canada.}
}
,
\, 
\textbf{J. C. Pardo
\footnote{Centro de Investigaci\'on en Matem\'aticas A.C. Calle Jalisco s/n. 36240 Guanajuato, M\'exico. Email: jcpardo@cimat.mx}
} 
}
\date{\footnotesize Current version: \today}

\begin{document}

\maketitle
\begin{abstract}
We study the distribution and various properties of exponential functionals of hypergeometric L\'evy processes.  We derive an explicit formula for the Mellin transform of the exponential functional and give both convergent and asymptotic series expansions of its probability density function. As applications we present a new proof of some of the results on the density of the supremum of a stable process, which 
were recently obtained in \cite{Kuznetsov2010} and \cite{KuzHub2010}. We also derive some new results related to (i) the entrance law of the stable
process conditioned to stay positive, (ii) the entrance law of the excursion measure of the stable process reflected at its past infimum and (iii) 
the entrance law and the last passage time of the radial part of n-dimensional symmetric stable process.

\bigskip


\medskip

\end{abstract}


\section{Introduction}
Exponential functionals of L\'evy
processes play a very important role in various domains of
probability theory, such as self-similar Markov processes, random processes in random environment, fragmentation processes, branching processes, mathematical finance, to name but a few. In general, the distribution of the exponential functional of a L\'evy process $X=(X_t,t\ge 0)$ with lifetime $\zeta$, defined as
\beqq
I:=\int\limits_0^{\zeta} e^{- X_t}\d t,
\eeqq
can be rather complicated. Nonetheless, this distribution is known explicitly for the case when $X$ is either: a standard Poisson processes,  Brownian motion with drift,  a particular class of spectrally negative Lamperti stable (see for instance \cite{CKP, KP,PP}), spectrally positive L\'evy processes satisfying the Cram\'er condition (see for instance \cite{PP1}). In the class of L\'evy processes with double-sided jumps the distibution of the exponential functional is known in closed form 
only in the case of L\'evy processes with hyperexponential jumps, see the recent paper by Cai and Kou \cite{CaiKou2010}. An overview of this topic can be found in Bertoin and Yor \cite{BYS}.

For many applications, it is enough to have estimates of $\p(I<t)$ as $t\to 0^+$ and $\p(I>t)$ as $t\to+\infty$. However it is quite difficult to 
obtain such estimates in the general case. The behaviour of the tail $\p(I>t)$  has been studied in a general setting (see for instance \cite{CKP,CP,HA,Maulik2006,RI,RI1}). On the other hand, the behaviour of $\p(I<t)$ has been studied only in two particular cases: when $X$ has exponential moments and its Laplace exponent is regularly varying at infinity with index $\theta\in (1,2)$ (see \cite{PA}) and  when $X$ is a subordinator whose Laplace exponent is regularly varying at $0$ (see \cite{CR}).

At the same time, the problem of finding the distribution of the supremum of a stable process has also stimulated a lot of research. The explcit expressions for the 
Wiener-Hopf factors for a dense set of 
parameters were first obtained by Doney \cite{Doney1987}. In the spectrally positive case a convergent series representation for the density of supremum 
was first obtained by Bernyk, Dalang and Peskir \cite{Bernyk2008}, and a complete asymptotic expantion was derived by Patie \cite{Patie2009}. The general case
was treated recently in \cite{Kuznetsov2010} and \cite{KuzHub2010}, where the authors have derived explicit formulas for the Wiener-Hopf factors, the Mellin transform of supremum and also convergent and asymptotic series representations for the density of the supremum. 

In this paper we pursue two goals. First, we will study exponential functionals of hypergeometric processes, we will obtain the Mellin transform and 
both convergent and asymptotic series representations for the density of the exponential functional. This gives us the first explicit results on the exponential functionals of L\'evy processes  which have double-sided jumps of infinite activity or infinite variation. The technique that we develop to obtain these results may be of independent interest, as it provides a very simple and straightforward method to derive well-known results on the exponential functionals 
of Brownian motion with drift and the less-known (but very interesting) recent results by Cai and Kou \cite{CaiKou2010} for processes with double-sided hyper-exponential jumps.

The second goal is to prove that hypergeometric processes include the Lamperti-stable processes, which will allow us to prove interesting results on
fluctuations of stable processes with the help of the Lamperti transformation. In particular, we will obtain a new proof of some of the results in 
 \cite{Kuznetsov2010} and \cite{KuzHub2010}, which is more straightforward and less technical. We will also derive 
some new formulas related to the density of the entrance law of the excursion measure of the 
stable process reflected at its past infimum and will obtain some several results related to the $n$-dimensional symmetric stable process.

The paper is organized as follows: in Section \ref{sec_hypergeometric} we introduce hypergeometric processes and establish the connections 
between this class and the Lamperti-stable processes. In Section \ref{sec_exponential_functionals} we study the Mellin transform of the exponential functional of the 
hypergeometric processes and in Section \ref{sec_density} we derive the convergent and asymptotic series representations for the density of the 
exponential functional. Finally, in Section \ref{sec_applications} we present some applications of these results to fluctuations of stable processes.

\section{Hypergeometric and Lamperti-stable processes}\label{sec_hypergeometric}

Hypergeometric processes were first introduced in \cite{KPR} and, more generally, in \cite{KKPS2010}. These processes were originally constructed using 
Vigon's theory of philantropy (see \cite{VI}) and they provide examples of L\'evy processes with an explicit Wiener-Hopf factorization. 
The class of processes which we will study in this paper should be considered as a generalization of the 
hypergeometric processes studied in \cite{KPR} and of the Lamperti-stable processes, which were introduced by Caballero and Chaumont \cite{CC}.

We start by defining a function $\psi(z)$ as
\begin{equation}\label{hype}
\psi(z)=- \frac{\Gamma(1-\beta+\gamma-z)}{\Gamma(1-\beta-z)}\frac{\Gamma(\hat\beta+\hat\gamma+z)}{\Gamma(\hat\beta+z)},
\end{equation}
where $(\beta,\gamma,\hat \beta,\hat \gamma)$ belong to the admissible set of parameters 
\beq\label{def_A}
{\mathcal A}=\{ \beta \le 1, \; \gamma \in (0,1), \; \hat \beta \ge 0, \; \hat \gamma \in (0,1) \}.
\eeq 
Our first goal is to show that $\psi(z)$ is the Laplace exponent of a (possibly killed) L\'evy process $X$, that is $\psi(z)=\ln \e [ \exp( z X_1) ]$. This process
and its properties will be presented in the next Proposition. From now on we will use the following notation
\beq\label{defn_eta}
 \eta=1-\beta+\gamma+\hat \beta+\hat \gamma.
\eeq

\begin{proposition}\label{prop1}

\noindent
\begin{itemize}
\item[(i)] The function $\psi(z)$ defined by (\ref{hype}) is the Laplace exponent of a L\'evy process $X$. The density of its L\'evy
   measure  is given by
 \beq\label{pi_formula}
 \pi(x)=
 \begin{cases}\displaystyle
-\frac{\Gamma(\eta)}{\Gamma(\eta-\hat \gamma)\Gamma(-\gamma)}
 e^{-(1-\beta+\gamma)x} {}_2F_1\left(1+\gamma,\eta;\eta-\hat \gamma;e^{-x}\right), \;\;\; {\textnormal{ if }} \; x>0, \\ {} \\
 \displaystyle 
 -\frac{\Gamma(\eta)}{\Gamma(\eta-\gamma)\Gamma(-\hat \gamma)}
 e^{(\hat\beta+\hat \gamma)x} {}_2F_1\left(1+\hat \gamma,\eta;\eta-\gamma;e^{x}\right), \;\;\;\;\qquad {\textnormal{ if }} \; x<0, 
 \end{cases}
 \eeq
where ${}_2F_1$ is the Gauss hypergeometric function.
\item[(ii)] When $\beta<1$ and $\hat \beta>0$ the process $X$ is killed at the rate 
\beqq
q=-\psi(0)=\frac{\Gamma(1-\beta+\gamma)}{\Gamma(1-\beta)}\frac{\Gamma(\hat\beta+\hat\gamma)}{\Gamma(\hat\beta)}.
\eeqq
When $\beta=1$ and $\hat \beta> 0$ \{$\beta<1$ and $\hat \beta=0$\}  $X$ drifts to $+\infty$ \{$-\infty$\} and
\beqq
\e[X_1]=\frac{\Gamma(\gamma)\Gamma(\hat \beta+\hat \gamma)}{\Gamma(\hat \beta)} \;\;\; \qquad
\Bigg\{ \e[X_1]=-\frac{\Gamma(\hat \gamma)\Gamma(1- \beta+\gamma)}{\Gamma(1-\beta)}\Bigg\}.
\eeqq
When $\beta=1$ and $\hat \beta=0$ the process $X$ oscillates.
\item[(iii)] The process $X$ has no Gaussian component. When $\gamma+\hat \gamma<1$ \{$1 \le \gamma+\hat \gamma<2$\} the process has paths of bounded variation and no linear drift
 \{paths of unbounded variation\}. 
\item[(iv)] We have the Wiener-Hopf factorization $-\psi(z)=\kappa(q,-z)\hat \kappa(q,z)$ where
\beq\label{def_kappa_kappahat}
\kappa(q,z)=\frac{\Gamma(1-\beta+\gamma+z)}{\Gamma(1-\beta+z)}, \;\;\; \hat \kappa(q,z)=\frac{\Gamma(\hat\beta+\hat\gamma+z)}{\Gamma(\hat\beta+z)}
\eeq
\item[(v)] The process $\hat X=-X$ is a hypergeometric process with parameters $(1-\hat \beta,\hat \gamma,1-\beta,\gamma)$. 
\end{itemize}
\end{proposition}
\begin{proof}
First let us prove (i).  Let $X^{(1)}$ be a general hypergeometric L\'evy process (see section 3.2 in \cite{KKPS2010}) with parameters
\beqq
&& \sigma=d=k_1=\delta_1=\delta_2=0, \; \beta=1, \; c_1=-\frac{1}{\Gamma(-\gamma)},  \; c_2=-\frac{1}{\Gamma(-\hat \gamma_2)}, \\
&& \; \alpha_1=\beta, \; \alpha_2=1-\hat \beta,\;  \gamma_1=\gamma,   \; \gamma_2=\hat \gamma, \;
 k_2=\frac{\Gamma(\hat \beta+\hat \gamma)}{\Gamma(\hat \beta)}.
\eeqq
This process is constructed using Vigon's theory of philantropy (see \cite{VI}) from two subordinators $H^{(1)}$ and $\hat H^{(1)}$ which have Laplace exponents
\beqq
\kappa^{(1)}(q,z)=\kappa(q,z)-\kappa(q,0), \;\;\; \hat \kappa^{(1)}(q,z)=\hat \kappa(q,z)
\eeqq
where $\kappa(q,z)$ and  $\hat \kappa(q,z)$ are given by (\ref{def_kappa_kappahat}). 
We see that the Lapalace exponent $\psi^{(1)}(z)$ of process $X^{(1)}$ satisfies
\beqq
 \psi^{(1)}(z)=\psi(z)+k \hat \kappa(q,z)
\eeqq
where we have denoted $k=\kappa(q,0)$. Therefore the process $X^{(1)}_t$ has the same distribution as $X_t-\hat H_{kt}$, in particular the distribution of 
positive jumps of $X^{(1)}$ coincides with the distribution of positive jumps of $X$. From \cite{KKPS2010} we find that the L\'evy measure of $X^{(1)}$ 
restricted to $x>0$ coincides with (\ref{pi_formula}). The expression of the L\'evy measure for $x<0$ follows easily by symmetry considerations. This proves that $\psi(z)$ defined by (\ref{hype}) is the Laplace transform of a (possibly killed) L\'evy process
with the density of the L\'evy measure given by (\ref{pi_formula}).

The rest of the proof is rather straightforward. Property (v) follows easily from the definition of the Laplace exponent (\ref{hype}). The Wiener-Hopf factorization (iv) follows easily by construction: we know that both $\kappa(q,z)$ and 
$\hat \kappa(q,z)$ defined by (\ref{def_kappa_kappahat}) are Laplace transforms of (possibly killed) subordinators, and the result follows 
from identity $-\psi(z)=\kappa(q,-z)\hat \kappa(q,z)$ and the uniqueness of the Wiener-Hopf factorization (see \cite{Kyprianou}). 

Let us prove (ii). The fact that $X$ drifts to $+\infty$ when $\beta=1$ and $\hat \beta>0$ follows from the Wiener-Hopf factorization (iv): in this case $\kappa(q,0)=0$, therefore the ascending ladder height process drifts to $+\infty$, while the descending ladder height process is killed at a rate $\hat \kappa(q,0)>0$. The expression
for $\e[X_1]$ follows from (\ref{hype}) using the fact that $\e[X_1]=\psi'(0)$. 
Other results in (ii) can be verified in a similar way. 

 Let us prove (iii). Formula 
(\ref{hype}) and the asymptotic expansion for the Gamma function (see formula 8.328.1 in \cite{Jeffrey2007}) imply that  
\beq\label{psi(iz)_asymptotics}
\psi(\i z)=O(|z|^{\gamma+\hat\gamma})=o(|z|^2), \;\;\; z\to \infty, z\;\in \r.
\eeq
 Applying Proposition 2 from \cite{Be} we conclude that $X$ has no Gaussian component, and that when $\gamma+\hat \gamma<1$ 
\{$1<\gamma+\hat \gamma<2$\} the process has paths of bounded variation and no linear drift \{paths of unbounded variation\}. In the remaining 
case $\gamma+\hat \gamma=1$ the density of the L\'evy measure has 
a singularity of the form $C x^{-2}+o(x^{-2})$ as $x\to 0^+$ (see formula 15.3.12 in \cite{AbramowitzStegun}), which implies that the process has 
paths of unbounded variation.
\end{proof}

{\bf Remark 1:} It is possible to extend the definition of the hypergeometric processes to the case when $\gamma \in \{0,1\}$ or $\hat \gamma \in \{0,1\}$. 
In the case when $\gamma=1$ and $\hat \gamma\in (0,1)$ the process $X$ is spectrally negative. This is true due to the formula (\ref{pi_formula}), which shows that 
$\pi(x)=0$ for $x>0$; one can also establish this fact using (\ref{def_kappa_kappahat}), which imlies that $\kappa(q,z)=1-\beta+z$. Similarly, when $\hat \gamma=1$ and 
$\gamma\in (0,1)$ the process $X$ is spectrally positive. When $\gamma=\hat\gamma=1$ the process $X$ is Brownian motion with drift, killed at an exponential time. When $\gamma=0$ \{$\hat \gamma=0\}$ the process $-X$ \{$X$\} is a subordinator. 
\\

The three Lamperti-stable processes $\xi^*,\xi^{\uparrow}$ and $\xi^{\downarrow}$ were introduced by Caballero and Chaumont \cite{CC} 
by applying the Lamperti transformation (see \cite{LA}) to the positive self-similar Markov processes related
to stable process. In particular, the process
$\xi^*$ is obtained from a stable process started at $x>0$ and killed on the first exit from the positive half-line, while the  process $\xi^{\uparrow}$ 
\{$\xi^{\downarrow}$\} is obtained from stable process conditioned to stay positive \{conditioned to hit zero continuously\}. We refer to \cite{CC, CKP, CPP} for all
the details on these processes.

For our next result we will need the form of the Laplace exponent of the Lamperti-stable process $\xi^*$. We assume that the characteristic exponent 
$\Psi_Y(z)=-\ln \e [ \exp(\i z Y_1) ]$ of a stable process $Y$ is defined as
\beq\label{def_Psi_stable}
\Psi_Y(z)=e^{\frac{\pi \i \gamma}2} z^{\alpha} {\bf 1}_{\{z>0\}}+e^{-\frac{\pi \i \gamma}2}|z|^{\alpha} {\bf 1}_{\{z<0\}},
\eeq
where $|\gamma| < \alpha$ if $\alpha \in (0,1)$ and $|\gamma|<1$ if $\alpha \in (1,2)$. 
With this parametrization the density of the L\'evy measure of $Y$ is given by
\beqq
\pi_Y(x) =   c_+  x^{-1-\alpha}{\mathbf 1}_{\{x>0\}} + 
 c_-  |x|^{-1-\alpha}{\mathbf 1}_{\{x<0\}},
\eeqq
where
\beq\label{def_cplus_cminus}
c_+=\Gamma(1+\alpha)\frac{\sin(\pi \alpha \rho)}{\pi}, \;\;\; c_-=\Gamma(1+\alpha)\frac{\sin(\pi \alpha(1-\rho))}{\pi}
\eeq
and $\rho=\p(Y_1>0)=(1- \gamma/\alpha)/2$. Note that this parameter was corresponds to  $1-\rho$ in \cite{CC}. According to Caballero and Chaumont \cite{CC}, the Laplace exponent $\psi^*(z)=\ln\e [ \exp(z \xi^*_1) ]$ of the 
Lamperti-stable process $\xi^*$ is given by
\begin{equation}\label{lexplam}
\psi^*(z)=\frac{c_+-c_-}{1-\alpha}z+\int_{\r\setminus\{0\}} \left(e^{z x}-1-z(e^x-1)\mathbf{1}_{\{|e^x-1|<1\}}\right)e^x\pi_Y\left(e^x-1 \right){\rm d}x -c_-\alpha^{-1}.
\end{equation} 
It is important to note that when $\alpha<1$, the Laplace exponent (\ref{lexplam}) can be rewritten as
\beqq
\psi^*(z)=\int_{\r\setminus\{0\}}  \left(e^{z x}-1\right)e^x\pi_Y\left(e^x-1 \right){\rm d}x -c_- \alpha^{-1},
\eeqq
so that in this case $\xi^*$ is a process of bounded variation with no linear drift.

\begin{theorem}\label{theorem_hyper=Lamperti_stable}
 Lamperti-stable processes $\xi^*$, $\xi^{\uparrow}$, $\xi^{\downarrow}$ can be identified as hypergeometric processes with the following sets of parameters 
\begin{table}[H]
\centering
\bigskip
\begin{tabular}{|c||c|c|c|c|}
\hline \rule[0pt]{0pt}{15pt}  
    & \;$\beta$ \;  & \; $\gamma$ \; & \; $\hat \beta$ \; & \; $\hat \gamma$ \;  \\ [0.5ex] \hline \hline 
   \rule[2pt]{0pt}{11pt}   $\xi^*$  & $1-\alpha (1-\rho)$ & $\alpha \rho$ & $1-\alpha (1-\rho)$ & $\alpha (1-\rho)$ \\ [0.5ex] \hline
   \rule[2pt]{0pt}{11pt}  $\xi^{\uparrow}$  & $1$ & $\alpha\rho$ & $1$ & $\alpha (1-\rho)$ \\ [0.5ex] \hline
   \rule[2pt]{0pt}{11pt}  $\xi^{\downarrow}$  & $0$ & $\alpha\rho$ & $0$ & $\alpha (1-\rho)$ \\ [0.5ex] \hline	
\end{tabular}
\end{table}
\end{theorem}
\begin{proof}
For the proof of the result for $\xi^{\uparrow}$ see Proposition 2 in \cite{KPR}. The result for $\xi^{\downarrow}$ follows from Proposition 1 in \cite{CKP}. Thus we only need to prove the result for $\xi^*$. 

Let us set $(\beta,\gamma,\hat \beta,\hat \gamma)=(1-\alpha (1-\rho), \alpha \rho, 1-\alpha (1-\rho), \alpha (1-\rho))$ and compute the L\'evy measure of the hypergeometric process $X$. We find that $\eta=1+\alpha$, and formula (\ref{pi_formula}) implies that for $x>0$ we have 
\beqq
 \pi(x)&=&-\frac{\Gamma(1+\alpha)}{\Gamma(1+\alpha \rho)\Gamma(-\alpha \rho)} e^{-\alpha x} {}_2F_1\left(1+\alpha \rho,1+\alpha;1+\alpha \rho;e^{-x}\right)\\
&=& \Gamma(1+\alpha)\frac{\sin(\pi\alpha \rho)}{\pi} e^{-\alpha x} (1-e^{-x})^{-1-\alpha}=e^x \pi_Y\left(e^x-1\right) 
\eeqq
where we have used the reflection formula for the gamma function and the fact that ${}_2F_1(a,b;a,z)=(1-z)^{-b}$.  Similarly, we find that
$\pi(x)=e^x \pi_Y(e^x-1)$ for $x<0$. We see that the L\'evy measure of the hypergeometric process is the same as the L\'evy measure of the Lamperti-stable
process. We know that in the case $\alpha<1$ both of these processes have paths of finite variation and no linear drift, this proves that $X\stackrel{d}{=}\xi^*$
in the case $\alpha<1$.

 In the case $\alpha>1$ processes $X$ and $\xi^*$ have infinite variation, no Gaussian component and identical L\'evy measures. 
Thus their Laplace exponents may differ only by a linear function. Thus in order to establish that the Laplace exponents are equal it is enough to show that
 $\e[X_1|\zeta>1]=\e[\xi^*_1|\zeta>1]$. Using (\ref{hype}) and the reflection formula for the gamma function we find that 
the Laplace exponent of the hypergeometric process $X$ is given by
\beqq
\psi_X(z)=\frac{1}{\pi}\Gamma(\alpha-z)\Gamma(1+z) \sin(\pi(z-\alpha(1-\rho)))
\eeqq 
Therefore we have
\beq\label{mean_X}
\e[X_1|\zeta>1]=\psi_X'(0)=\Gamma(\alpha) \frac{\sin(\pi \alpha (1-\rho))}{\pi} (\Psi(\alpha)-\Psi(1))+\Gamma(\alpha)\cos(\pi \alpha  (1-\rho)).
\eeq
where $\Psi(z)=\Gamma'(z)/\Gamma(z)$ is the digamma function. On the other hand, equation (11) in \cite{CC} tells us that 
\beq\label{mean_xi*}
\e[\xi^*_1|\zeta>1]=\frac{c_+-c_-}{1-\alpha}&+&c_+\left[ \int\limits_0^{\ln(2)} \frac{(1+y-e^y)e^y}{(e^y-1)^{1+\alpha}} \d y+
 \int\limits_{\ln(2)}^{\infty} \frac{ye^y}{(e^y-1)^{1+\alpha}} \d y  \right] \\ \nonumber
&+&c_-\int\limits_{-\infty}^0 \frac{(1+y-e^y)e^y}{(1-e^y)^{1+\alpha}}\d y,
\eeq
where $c_+$ and $c_-$ are defined in (\ref{def_cplus_cminus}). Let us show that the expressions in the right-hand side of (\ref{mean_X}) and (\ref{mean_xi*}) are equal.

First we will deal with the integrals multiplying $c^+$ in (\ref{mean_xi*}). We rearrange the terms as follows
\beq\label{lam_stable_proof1}
\int\limits_0^{\ln(2)} \frac{(1+y-e^y)e^y}{(e^y-1)^{1+\alpha}} \d y+
 \int\limits_{\ln(2)}^{\infty} \frac{ye^y}{(e^y-1)^{1+\alpha}} \d y =
\int\limits_0^{\infty} \frac{(1+y-e^y)e^y}{(e^y-1)^{1+\alpha}} \d y+
 \int\limits_{\ln(2)}^{\infty} \frac{(e^y-1)e^y}{(e^y-1)^{1+\alpha}} \d y
\eeq
Performing the change of variables $u=\exp(y)$ it is easy to see that the second integral in the right-hand side of (\ref{lam_stable_proof1}) is equal to $1/(\alpha-1)$. In order to compute the first integral we use integration by parts
\beq\label{lam_stable_proof2}
\int\limits_0^{\infty} \frac{(1+y-e^y)e^y}{(e^y-1)^{1+\alpha}} \d y=
-\frac{1}{\alpha} \int\limits_0^{\infty} (1+y-e^y) \d (e^y-1)^{-\alpha} =
\frac{1}{\alpha} \int\limits_0^{\infty} \frac{1-e^y}{(e^y-1)^{\alpha}}\d y=\frac{\pi}{\alpha\sin(\pi \alpha)}
\eeq
where in the last step we have used the following integral formula
\beq\label{beta_integal}
\int\limits_0^{\infty} \frac{e^{a u}}{(e^{b u}-1)^{c}} \d u=\frac{1}{b} \frac{\Gamma\left(c-\frac{a}{b}\right)\Gamma(1-c)}{\Gamma\left(1-\frac{a}{b}\right)}, 
\;\;\; \frac{a}{b}<c<1,
\eeq
and the reflection formula for the gamma function.
The last integral in (\ref{mean_xi*}) can also be computed using integration by parts
\beq\label{lam_stable_proof3}
&&\int\limits_{-\infty}^0 \frac{(1+y-e^y)e^y}{(1-e^y)^{1+\alpha}}\d y=\int\limits_0^{\infty} \frac{(1-y-e^{-y})e^{\alpha y}}{(e^y-1)^{1+\alpha}}\d y\\
\nonumber
&& \qquad =
-\frac{1}{\alpha} \int\limits_0^{\infty} (1-y-e^{-y})e^{(\alpha-1)y} \d (e^y-1)^{-\alpha}\\ \nonumber
&& \qquad =
\frac{1}{\alpha} \int\limits_0^{\infty} \frac{(-1+e^{-y})e^{(\alpha-1)y}+(1-y-e^{-y})(\alpha-1)e^{(\alpha-1)y} }{(e^y-1)^{\alpha} }\d y\\ \nonumber
&& \qquad =
\frac{\alpha-2}{\alpha}\int\limits_{0}^{\infty} \frac{e^{(\alpha-2)y}}{(e^y-1)^{\alpha-1}}\d y - 
\frac{\alpha-1}{\alpha} \int\limits_0^{\infty} \frac{ye^{(\alpha-1)y}}{(e^y-1)^{\alpha}}\d y\\ \nonumber
&& \qquad =
 \frac{\alpha-2}{\alpha} {\textnormal{B}}(2-\alpha,1)- \frac{\alpha-1}{\alpha} \frac{\d}{\d z} \left[  \int\limits_0^{\infty} \frac{e^{z y}}{(e^y-1)^{\alpha}}\d y \right]_{z=\alpha-1}\\ \nonumber
&& \qquad = -\frac{1}{\alpha}-\frac{1}{\alpha}(\Psi(1)-\Psi(2-\alpha))
\eeq
where in the last step we have again used (\ref{beta_integal}). Combining (\ref{mean_xi*}), (\ref{lam_stable_proof1}), (\ref{lam_stable_proof2}) and (\ref{lam_stable_proof3}) we see that
\beqq
 \e[\xi^*_1|\zeta>1]=\frac{c_+-c_-}{1-\alpha}+c_+\left[ \frac{1}{\alpha-1}+\frac{\pi}{\alpha\sin(\pi\alpha)} \right]+
c_-\left[-\frac{1}{\alpha}-\frac{1}{\alpha}(\Psi(1)-\Psi(2-\alpha))\right].
\eeqq
Using the fact that $\Psi(2-\alpha)=\Psi(\alpha)+\pi\cot(\pi\alpha)+1/(1-\alpha)$, the identity $\Gamma(z+1)=z\Gamma(z)$ and addition formula for the sine function it is not hard to reduce the above expression to (\ref{mean_X}), we leave all the details to the reader.
\end{proof}

\section{Mellin transform of the exponential functional}\label{sec_exponential_functionals}

Let $X$ be the hypergeometric L\'evy process with parameters $(\beta,\gamma,\hat\beta,\hat \gamma) \in {\mathcal A}$ and $\hat \beta>0$.  
We assume that $\alpha>0$ and define the exponential functional as
\beq\label{def_exp_functional}
 I(\alpha,X)=\int\limits_{0}^{\zeta} e^{-\alpha X_t} \d t,
\eeq
where $\zeta$ is the lifetime of the process $X$. Note that the above integral converges with probability one: according to Proposition \ref{prop1}, either $\zeta$ is finite (if $\beta<1$) or $\zeta=+\infty$ and the process $X$
drifts to $+\infty$ (if $\beta=1$).

Our main tool in studying the exponential functional is the Mellin transform, which is defined as
\beq\label{def_Mellin_transform}
\mm(s)=\mm(s;\alpha,\beta,\gamma,\hat \beta,\hat \gamma)=\e[I(\alpha,X)^{s-1}].
\eeq
From the definition of the Laplace exponent (\ref{hype}) we find that $X$ satisfies Cram\'er's condition $\e[\exp(-\hat \beta X_1)]=1$, therefore applying 
Lemma 2 from \cite{r2007} we conclude that $\mm(s)$ exists for $s \in (0,1+\hat \beta \delta)$. 

In order to describe our main result in this section, we need to define the Barnes double gamma function, $G(z;\tau)$. This function was introduced by 
Alexeiewsky in 1889 and was extensively studied by Barnes in \cite{Barnes1899} and \cite{Barnes1901}. The double gamma function is defined  as an infinite product in Weierstrass's form 
\beq\label{def_G_z_tau}
G(z;\tau)=\frac{z}{\tau} e^{a\frac{z}{\tau}+b\frac{z^2}{2\tau}} \prod\limits_{m\ge 0} \prod\limits_{n\ge 0} {}^{'}
\left(1+\frac{z}{m\tau+n} \right)e^{-\frac{z}{m\tau+n}+\frac{z^2}{2(m\tau+n)^2}}, \;\;\; 
 |\arg(\tau)|<\pi, \;\;\; z\in \c. 
\eeq
Here the prime in the second product means that the term corresponding to $m=n=0$ is omitted.
Note that by definition $G(z;\tau)$ is an entire function in $z$ and has simple zeros on the lattice $m\tau+n$, $m\le 0$, $n\le 0$. 
The function $G(z;\tau)$ can also be expressed as an infinite product  of gamma functions, see \cite{Barnes1899}. An integral representation for $\ln(G(z;\tau))$ was obtained in \cite{Lawrie1994} and several important asymptotic expansions were established in \cite{Bill1997}. 

It turns out that it is possible to define constants $a=a(\tau)$ and $b=b(\tau)$ in a particular way, so that we have $G(1;\tau)=1$ and that $G(z;\tau)$ satisfies the following three functional identities (see \cite{Barnes1899} and \cite{Kuznetsov2010})
\beq\label{funct_rel_G_1}
G(z+1;\tau)&=&\Gamma\left(\frac{z}{\tau}\right) G(z;\tau), \\
\label{funct_rel_G_tau}
 G(z+\tau;\tau)&=&(2\pi)^{\frac{\tau-1}2}\tau^{-z+\frac12}
\Gamma(z) G(z;\tau),\\
\label{eq_G_1_over_tau}
G(z;\tau)&=&(2\pi)^{\frac{z}2 \left(1-\frac1{\tau} \right)} \tau^{-\frac{z^2}{2\tau}+\frac{z}{2} \left( 1+ \frac{1}{\tau} \right)-1 }  G\left(\txt\frac{z}{\tau};\txt\frac{1}{\tau}\right).
\eeq

Next we introduce a function which will play a central role in our study of exponential functionals. Everywhere in this paper we will use notation
\beqq
\delta=\frac{1}{\alpha}.
\eeqq
\begin{definition}\label{def_Ms} For $s\in \c$ we define
 \beq\label{def_Mellin_scaled}
M(s)=M(s;\alpha,\beta,\gamma,\hat \beta,\hat \gamma)=C \frac{G ((1-\beta)\delta+s; \delta)}{G ((1-\beta+\gamma)\delta+s;\delta) }
 \frac{G((\hat \beta+\hat \gamma)\delta+1-s;\delta)}{G(\hat \beta \delta +1-s;\delta)}, \;\;\; s\in \c.
\eeq
where the constant $C$ is such that $M(1)=1$.
\end{definition}

Our main result in this section is the following Theorem, which provides an explicit expression for the Mellin transform of the exponential functional.
\begin{theorem}\label{thm_main} 
Assume that $\alpha>0$, $(\beta,\gamma,\hat\beta,\hat \gamma) \in {\mathcal A}$ and $\hat \beta>0$. Then $\mm(s)\equiv \Gamma(s)M(s)$ for all $s\in \c$.
\end{theorem}

Before we are able to prove Theorem \ref{thm_main},  we need to establish several auxiliary results.
\begin{lemma}\label{lemma_G2_asymptotics}

\noindent
\begin{itemize}
 \item[(i)]  Assume that $\tau>0$. When $s\to \infty$ in the domain $|\arg(s)|<\pi-\epsilon<\pi$ we have
 \beq\label{eqn_G2_asymptotics}
 \ln \left[ \frac{G(a+s;\tau)}{G(s;\tau)} \right]= \frac{a}{\tau} s \ln(s) - \frac{a}{\tau} (1+\ln(\tau)) s + \frac{a}{2\tau} (a-1-\tau) \ln(s)+O(1).
 \eeq
 \item[(iv)] When $s\to \infty$ in the domain $0<\epsilon<\arg(s)<\pi - \epsilon$ we have
 \beq\label{eqn_Ms_asymptotics}
  \ln(M(s))=-(\gamma+\hat\gamma) s \ln(s) + s ((1+\ln(\delta))(\gamma+\hat \gamma)+\pi \i \hat \gamma)+O(\ln(s)).
 \eeq
 \item[(iii)] When $s\to \infty$ in a vertical strip $a < \re(s) < b$  we have
 \beq\label{eqn_Ms_asymptotics_Ims}
 \big | M(s) \big |=\exp\left( \frac{\pi}2 (\gamma-\hat \gamma) |\im(s)|+O(\ln|\im(s)|) \right).
\eeq
\end{itemize}
\end{lemma}
\begin{proof}
 Part (i) follows from the asymptotic expansion for $G(z;\tau)$, given in formula (4.5) in \cite{Bill1997}, while parts (ii) and (iii) are simple corollaries of 
 (i) and (\ref{def_Mellin_scaled}). 
\end{proof}

The following notation will be used extensively in this paper: if $X$ is a hypergeometric process with parameters $(\beta,\gamma,\hat \beta,\hat \gamma)$, then
 $\tilde X$ denotes the hypergeometric process with parameters $(\delta \beta, \delta \gamma,\delta\hat\beta, \delta\hat \gamma)$, 
provided that this parameter set is admissible. In particular, the Laplace exponent of $\tilde X$ is given by
\beq\label{def_tilde_psi}
\tilde\psi(z)=
-\frac{\Gamma(1-\delta(\beta-\gamma)-z) \Gamma(\delta(\hat \beta + \hat \gamma)+z)}
{\Gamma(1-\delta\beta-z)  \Gamma(\delta\hat \beta+z)}.
\eeq

\begin{lemma}\label{lemma_Ms_properties}
 
\noindent
\begin{itemize}
 \item[(i)] $M(s)$ is a real meromorphic function which has zeros 
\beq\label{M_s_zeros}
\bigg \{-(1-\beta)\delta-m\delta - n, \;\;\; 1+(\hat \beta +\hat \gamma)\delta+m\delta+n \bigg \}_{m,n\ge 0},
\eeq
and poles
\beq\label{M_s_poles}
\bigg \{z^{-}_{m,n}=-(1-\beta+\gamma)\delta-m\delta-n, \;\;\; z^{+}_{m,n}=1+\hat \beta\delta+m\delta + n  \bigg \}_{m,n\ge 0}. 
\eeq
 All zeros/poles are simple if $\alpha \notin {\mathbb Q}$. 
\item[(ii)]
$M(s)$ satisfies the following functional identities
\beq\label{eqn_M_s_plus_1}
M(s+1)&=&-\frac{1}{\psi(-\alpha s)}M(s), \\
 \label{eqn_M_s_plus_delta}
M(s+\delta)&=&- \frac{\alpha^{-\delta(\hat \gamma+\gamma)}}{\tilde\psi(1-\delta-s)}M(s),\\
\label{eqn_M_transform}
M(s;\alpha,\beta,\gamma,\hat \beta,\hat \gamma)&=&
\alpha^{(1-s)(\gamma+ \hat\gamma)} M(1-\alpha+\alpha s;\delta,\delta \beta,\delta \gamma,\delta \hat \beta,\delta \hat \gamma).
\eeq
\end{itemize}
\end{lemma}
\begin{proof}
 The proof of (i) follows from the definition of the double gamma function (\ref{def_G_z_tau}), while functional identity (\ref{eqn_M_s_plus_1}) 
is a simple corollary of the functional identity for the double gamma function (\ref{funct_rel_G_1}).  

Let us prove (\ref{eqn_M_transform}). We use (\ref{eq_G_1_over_tau}) and find that for all $s$ and $x$ 
\beqq
\frac{G\left( s+x ;\delta \right)}{G\left(s;\delta \right)}
= \alpha^{\alpha s x} C \frac{G(\alpha(s+x);\alpha)}{G(\alpha s;\alpha)}
\eeqq
where $C=C(x)$ depends only on $x$. The above identity implies that  
\beqq
&&\frac{G ((1-\beta)\delta+s; \delta)}{G ((1-\beta+\gamma)\delta+s;\delta) }
 \frac{G((\hat \beta+\hat \gamma)\delta+1-s;\delta)}{G(\hat \beta \delta +1-s;\delta)}\\
&&\qquad\qquad=\tilde C \alpha^{-s(\gamma+\hat\gamma)} \frac{G (1-\beta + \alpha s; \alpha)}{G (1-\beta+\gamma+ \alpha s;\alpha) }
 \frac{G(\hat \beta+\hat \gamma+\alpha-\alpha s;\alpha)}{G(\hat \beta  +\alpha - \alpha s;\alpha)}
\eeqq
where $\tilde C$ does not depend on $s$. It is easy to see that (\ref{eqn_M_transform}) follows from the above identity and (\ref{def_Mellin_scaled}).

The functional identity (\ref{eqn_M_s_plus_delta}) follows from (\ref{eqn_M_transform}) and (\ref{eqn_M_s_plus_1}), we leave all the details to the reader.
\end{proof}

The following proposition will be central in the proof of Theorem \ref{thm_main}. It allows us to identify explicitly the Mellin transform of the exponential
functional. This result is also applicable to some other L\'evy process, including Brownian motion with drift and more generally, processes with hyperexponential or phase-type jumps, therefore it is of independent interest. 

 First let us present the main ingredients.  Let $Y$ be a (possibly killed) L\'evy process started from zero, and let $\psi_Y(z)=\ln\e[\exp(z Y_1)]$ denote its Laplace exponent. In the case when $\psi_Y(0)=0$ (the process is not killed) we will also assume that $\e[Y_1]>0$, so that $Y$ drifts to $+\infty$. As usual we 
define the exponential functional $I=\int_0^{\zeta} \exp(-Y_t) \d t$ (where $\zeta$ is the lifetime of $Y$) and the Mellin transform $\mm_Y(s)=\e[I^{s-1}]$. 

\begin{proposition}\label{Ms_uniqueness_lemma} {\bf (Verification result)}
Assume that Cram\'er's condition is satisfied: there exists $z_0<0$ such that $\psi_Y(z)$ is finite for all $z\in (z_0, 0)$ and $\psi_Y(-\theta)=0$ for some $\theta \in (0,-z_0)$. If $f(s)$ satisfies the following three properties
\begin{itemize}
 \item[(i)] $f(s)$ is analytic and zero-free in the strip $\re(s)\in (0,1+\theta)$,
 \item[(ii)] $f(1)=1$ and $f(s+1)=-s f(s)/\psi_Y(-s)$ for all $s\in (0,\theta)$,  
 \item[(iii)] $|f(s)|^{-1}=o(\exp(2 \pi |\im(s)|))$ as $\im(s)\to \infty$, $\re(s)\in (0,1+\theta)$,
\end{itemize}
 then $\mm_Y(s) \equiv f(s)$ for  $\re(s) \in (0,1+\theta)$.
\end{proposition}
 \begin{proof}
The Cram\'er's condition and Lemma 2 in \cite{r2007} imply that $\mm_Y(s)$ can be extended to an analytic function in the strip $\re(s) \in (0,1+\theta)$. 
In the case when $\psi_Y(0)=0$ \{$\psi_Y(0)<0$\} we use Lemma 2.1 in \cite{Maulik2006} \{Proposition 3.1 from \cite{CarPetYor1997}\} to conclude that $\mm_Y(s)$ satisfies the functional identity 
\beq\label{M_Y_funct_eqn}
\mm_Y(s+1)=-\frac{s}{\psi_Y(-s)} \mm_Y(s)
\eeq
 for all $s\in (0,\theta)$. Since $f(s)$ satisfies the same functional identity we conclude that the function
$F(s)=\mm_Y(s)/f(s)$ satisfies
$F(s+1)=F(s)$ for all $s\in (0,\theta)$. 
Using the assumption that $f(s)$ is analytic and zero-free we conclude that
$F(s)$ is an analytic function in the strip $\re(s) \in (0,1+\theta)$. Since $F(s)$ is also periodic with period equal to one, it can be extended to an analytic and periodic function in the entire complex plane. 

Our goal now is to prove that function $F(s)$ is constant. Since $F(s)$ is analytic and periodic in the entire complex plane, it can be represented as a Fourier series
\beq\label{g_fourier_series}
F(s)=c_0+F_1(\exp(2\pi \i s))+F_2(\exp(-2\pi i s)),
\eeq 
where $F_i(z)$ are entire functions
\beqq
F_1(z)=\sum\limits_{n\ge 1} c_n z^n, \;\;\; F_2(z)=\sum\limits_{n\ge 1} c_{-n} z^n
\eeqq
 Due to the inequality $|\mm_Y(s)|<\mm_Y(\re(s))$ and assumption (iii) we conclude that
$F(s)=o(\exp(2 \pi |\im(s)|))$ as $\im(s)\to \infty$. As $\im(s) \to -\infty$ this estimate and equation (\ref{g_fourier_series}) imply that 
 $F_1(z)=o(z)$ as $z\to \infty$, and using Cauchy's estimates 
(Proposition 2.14 in \cite{Conway1978}) we conclude that $F_1(z)\equiv 0$. Similarly we find that $F_2(z)\equiv 0$, thus $F(s)$ is a constant, 
and the  value of this constant is one, since $\mm_Y(1)=f(1)=1$.
 \end{proof}

We would like to stress that Proposition \ref{Ms_uniqueness_lemma} is an important result of indepent interest. We know that if Cram\'er's condition is satisfied then the Mellin transform $\mm_Y(s)$ satisfies the functional identity (\ref{M_Y_funct_eqn}), however it is clear that there are infinitely many functions 
which satisfy the same functional identity. Proposition \ref{Ms_uniqueness_lemma} tells us that 
if we have found such a function $f(s)$, which satisfies (\ref{M_Y_funct_eqn}), and if we can verify the two conditions about the zeros of this function
and its asymptotic behaviour, then we can in fact {\it uniquely  identify} $\mm_Y(s)\equiv f(s)$. In particular, this proposition can be used to provide a very 
simple and short proof of the well-known result on exponential functional of Brownian motion with drift and of the recent results on
exponential functionals of processes with double-sided hyper-exponential jumps (see \cite{CaiKou2010}).
\\

{\it Proof of Theorem \ref{thm_main}:}
First of all, we check that the Cram\'er's condition is satisfied with $\theta=\hat \beta \delta$. 
Let $f(s)=\Gamma(s)M(s)$, where $M(s)$ is defined by (\ref{def_Mellin_scaled}). From Lemma
\ref{lemma_Ms_properties}(i) we know that $f(s)$ is analytic and zero-free in the strip $\re(s) \in (0,1+\hat \beta\delta)$. By construction we have $f(1)=1$, and 
from the formula (\ref{eqn_M_s_plus_1}) we find that $f(s)$ satisfies $f(s+1)=-sf(s)/\psi(-\alpha s)$ for $s\in (0,\hat \beta\delta)$. Next, Lemma
\ref{lemma_G2_asymptotics}(iii) and Stirling's asymptotic formula for the gamma function (see formula 8.327.3 in \cite{Jeffrey2007}) 
\beq\label{stirling}
\ln(\Gamma(s))=s\ln(s)-s+O(\ln(s))
\eeq 
imply that as 
$s\to \infty$ in the vertical strip $\re(s)\in (0,1+\hat \beta\delta)$ we have 
\beqq
|f(s)|^{-1}=\exp\left(\frac{\pi}2 \left(1-\gamma+\hat \gamma \right) |\im(s)|+o(\im(s)) \right)=o\left(
 \exp(\pi |\im(s)|) \right),
\eeqq 
where in the last step we have also used the fact that both $\gamma$ and $\hat \gamma$ belong to the interval $(0,1)$.

We see that function $f(s)$ satisfies all conditions of Proposition \ref{Ms_uniqueness_lemma}, thus we can conclude that $\mm(s)\equiv f(s)$.
\qed

\begin{corollary}
Assume that $\alpha>0$, $\hat \beta>0$ and that both sets of parameters $(\beta, \gamma, \hat \beta,\hat \gamma)$ and 
$( \delta\beta, \delta \gamma, \delta \hat \beta, \delta \hat \gamma)$ belong to the admissible set ${\mathcal A}$. 
Then we have the following identity in distribution
\beq\label{exp_func_dist_identity}
\epsilon_1^{\alpha} \times I(\alpha;X ) \stackrel{d}{=}
 \alpha^{\gamma +\hat \gamma} \times \epsilon_1 \times I(\delta; \tilde X)^{\alpha}
\eeq
where $\epsilon_1 \sim Exp(1)$ and all random variables are assumed to be independent.
\end{corollary}
\begin{proof}
Rewrite (\ref{eqn_M_transform}) as
\beqq
\Gamma(1-\alpha+\alpha s)\mm(s;\alpha,\beta,\gamma,\hat \beta,\hat \gamma)=
\alpha^{(s-1)(\gamma+\hat\gamma)} \Gamma(s)\mm(1-\alpha+\alpha s;\delta,\delta \beta,\delta \gamma,\delta \hat \beta,\delta \hat \gamma)
\eeqq
and use the following facts: (i) $\Gamma(s)=\e [\epsilon_1^{s-1} ]$; (ii) if $f(s)$ is the Mellin transform of a random variable $\xi$ then
$f(1-\alpha+\alpha s)$ is the Mellin transform of the random variable $\xi^{\alpha}$; (iii) Mellin transform of the product of independent random variables
is the product of their Mellin transforms.
\end{proof}

\section{Density of the exponential functional}\label{sec_density}

In this section we will study the density of the exponential functional, defined as
\beqq
 p(x)=\frac{\d }{\d x} \p(I(\alpha,X) \le x), \;\;\; x\ge 0.
\eeqq
As in the previous section, $X$ is a hypergeometric process with parameters $(\beta,\gamma,\hat\beta,\hat \gamma) \in {\mathcal A}$ and $\hat \beta>0$. The main results of this section are the convergent series representations and complete asymptotic expansions of this function as $x \to 0^+$ or $x\to +\infty$. 

Let us define the following three sets of parameters which will be used extensively later. In the following definition (and everywhere else in this paper) $\psi(\cdot)$, $\tilde \psi(\cdot)$ and $M(\cdot)$ denote the functions which were defined in (\ref{hype}), (\ref{def_tilde_psi}) and (\ref{def_Mellin_scaled}); the sequences $z^-_{m,n}$ and $z^+_{m,n}$ represent the poles of $M(\cdot)$ and
were defined in (\ref{M_s_poles}); the constant $\eta$ is defined by (\ref{defn_eta}).  
\begin{definition}\label{def_an_bmn_cmn}
 Define the coefficients $\{a_n\}_{n\ge 0}$ as
\beq\label{def_an}
a_{n}=-\frac{1}{n!} \prod\limits_{j=0}^n \psi(\alpha j), \;\;\; n\ge 0.
\eeq
The coefficients $\{b_{m,n}\}_{m,n\ge 0}$ are defined recursively
 \beq\label{def_bmn}
 \begin{cases}
 b_{0,0}& = \displaystyle \delta\frac{\Gamma(\eta)\Gamma(-(1-\beta+\gamma)\delta) }{\Gamma(\eta-\hat \gamma) \Gamma(- \gamma)} M(1-(1-\beta+\gamma)\delta),
 \vspace{0.2cm}
\\ 
 b_{m,n}&=\displaystyle -\frac{\psi(-\alpha z^{-}_{m,n}) }{z^-_{m,n}}  b_{m,n-1}, \;\;\; m\ge0, \; n\ge 1, 
 \vspace{0.2cm}
\\
 b_{m,n}&=\displaystyle -\alpha^{\delta(\gamma+\hat \gamma)} \tilde \psi(1-\delta-z^{-}_{m,n}) \frac{\Gamma(z^{-}_{m,n})}{\Gamma(z^{-}_{m-1,n})}b_{m-1,n} , \;\;\; m\ge 1, \; n\ge 0.
 \end{cases}
 \eeq
Similarly, $\{c_{m,n}\}_{m,n\ge 0}$ are defined recursively
 \beq\label{def_cmn}
 \begin{cases}
 c_{0,0}&=\displaystyle\delta \frac{\Gamma(1+\hat \beta \delta) \Gamma(1-\beta+\hat \beta)}{\Gamma(\eta-\hat \gamma) \Gamma(\hat \gamma)}
 M(\hat \beta \delta)
\vspace{0.2cm}
\\
c_{m,n}&=\displaystyle - \frac{z^{+}_{m,n-1}}{ \psi(-\alpha z^{+}_{m,n-1})} c_{m,n-1}, \;\;\;, m\ge0, \; n\ge 1,
\vspace{0.2cm}
\\
c_{m,n}&=\displaystyle - \frac{\alpha^{-\delta(\gamma+\hat \gamma)}}{\tilde \psi(1-\delta-z^{+}_{m-1,n})} \frac{\Gamma(z^{+}_{m,n})}{\Gamma(z^{+}_{m-1,n})} c_{m-1,n}, \;\;\; m\ge 1, \;n\ge 0.
 \end{cases}
 \eeq
\end{definition}

Note that if $\beta=1$ we have $\psi(0)=0$, which implies that $a_n=0$ for all $n \ge 0$.  

\begin{proposition}\label{prop_residues} Assume that $\alpha \notin {\mathbb Q}$. For all $m,n\ge 0$ we have
 \beqq
 {\textnormal{Res}}(\mm(s) :  s=-n)&=&a_n, \;\;\; \textnormal{if} \; \beta<1, \\
  {\textnormal{Res}}(\mm(s) :  s=z^-_{m,n})&=&b_{m,n}, \\
 {\textnormal{Res}}(\mm(s) :  s=z^+_{m,n})&=&-c_{m,n}. 
 \eeqq
\end{proposition}
\begin{proof}
 Let us prove that the residue of $\mm(s)$ at $s=z^{-}_{m,n}$ is equal to $b_{m,n}$. First, 
we use Theorem \ref{thm_main} and rearrange the terms in the  functional identity (\ref{eqn_M_s_plus_1}) to find that
\beqq
\mm(s)=\frac{\delta M(s+1)}{s+(1-\beta+\gamma)\delta}  \frac{\Gamma(2-\beta+\gamma+\alpha s)\Gamma(\hat \beta+\hat \gamma-\alpha s)}
{\Gamma(1-\beta+\alpha s)\Gamma(\hat \beta-\alpha s)}.
\eeqq
The above identity and definition (\ref{def_bmn}) imply that as $s\to -(1-\beta+\gamma)\delta$
\beqq
\mm(s)=\frac{b_{0,0}}{s+(1-\beta+\gamma)\delta}+O(1),
\eeqq
which means that the residue of $\mm(s)$ at $z^-_{0,0}=-(1-\beta+\gamma)\delta$ is equal to $b_{0,0}$. 

Next, let us prove that the residues satisfy the second recursive identity in (\ref{def_bmn}). We rewrite (\ref{eqn_M_s_plus_1}) as
\beq\label{proof_residues}
 \mm(s)=-\frac{\psi(-\alpha s)}{s} \mm(s+1). 
\eeq
We know that $\mm(s)$ has a simple pole at $s=z^-_{m,n}$  while $\mm(s+1)$ has a simple pole at $z^-_{m,n}+1=z^-_{m-1,n}$. One can also check that
function $\psi(-\alpha s)$ is analytic at 
$s=z^-_{m,n}$ for $m\ge 1$. Therefore we have as $s\to z^-_{m,n}$
\beqq
\mm(s)&=&{\textnormal{Res}}(\mm(s):s=z^-_{m,n})\frac{1}{s-z^-_{m,n}}+O(1), \\
\mm(s+1)&=& {\textnormal{Res}}(\mm(s):s=z^-_{m-1,n})\frac{1}{s-z^-_{m,n}}+O(1),\\
-\frac{\psi(-\alpha s)}{s}&=&-\frac{\psi(-\alpha z^-_{m,n})}{z^-_{m,n}}+O(s-z^-_{m,n})
\eeqq
which, together with (\ref{proof_residues}) imply that
\beqq
{\textnormal{Res}}(\mm(s):s=z^-_{m,n})=-\frac{\psi(-\alpha z^-_{m,n})}{z^-_{m,n}} \times {\textnormal{Res}}(\mm(s):s=z^-_{m-1,n}).
\eeqq

The proof of all remaining cases is very similar and we leave the details to the reader.
\end{proof}

Proposition \ref{prop_residues} immediately gives us a complete asymptotic expansion of $p(x)$ as $x\to 0^+$ and $x\to +\infty$, presented in the next Theorem. 
\begin{theorem}\label{thm_density_asymptotics} Assume that $\alpha \notin {\mathbb Q}$. Then
 \beq\label{kx_asympt_x_to_0}
 p(x) &\sim& \sum\limits_{n\ge 0} a_n x^n + \sum\limits_{m\ge 0} \sum\limits_{n\ge 0} b_{m,n} x^{(m+1-\beta+\gamma)\delta+n}, \;\;\; x\to 0^+, \\
 \label{kx_asympt_x_to_infty}
 p(x) &\sim& \sum\limits_{m\ge 0} \sum\limits_{n\ge 0} c_{m,n} x^{-(m+\hat \beta)\delta-n-1}, \qquad\qquad\;\;\;\;\;\; x\to +\infty.
 \eeq
\end{theorem}
\begin{proof}

  The starting point of the proof is the expression of $p(x)$ as the inverse Mellin transform
\beq\label{thm4_proof0}
 p(x) = \frac{1}{2\pi \i} \int\limits_{1+\i \r} \mm(s) x^{-s} \d s, \;\;\; x>0.
\eeq
Due to (\ref{eqn_Ms_asymptotics_Ims}), Theorem \ref{thm_main} and Stirling's formula (\ref{stirling}) we know that $|\mm(x+\i u)|$ decreases exponentially as $u\to \infty$ (uniformly in $x$ in any finite interval), therefore the integral in the right-hand side of (\ref{thm4_proof0}) converges absolutely and $p(x)$ is a smooth function for $x>0$. 
Assume that $c<0$ satisfies $c \ne z^-_{m,n}$ and $c \ne -n$ for all $m,n$. Shifting the contour of integration $1+\i \r \mapsto c + \i \r$ and taking into account the residues at the poles $s=z^{-}_{m,n}$, we find that 
\beq\label{p_shift_contour}
p(x)=\sum \textnormal{Res}(\mm(s): s=z^-_{m,n}) \times x^{-z^-_{m,n}}&+&
\sum\limits_{0\le n < |c|} \textnormal{Res}(\mm(s): s=-n) \times x^{n} \\ \nonumber
&+& \frac{1}{2\pi \i} \int\limits_{c+\i \r} \mm(s) x^{-s} \d s.
\eeq
where the first summation is over all $m\ge 0, n \ge 0$, such that $z^-_{m,n}>c$. Next, we perform a change of variables $s=c+\i u$ and obtain the following estimate
\beqq
\Bigg |\int\limits_{c+\i \r} \mm(s) x^{-s} \d s \Bigg | = x^{-c} \Bigg |\int\limits_{\r} \mm(s) x^{-\i u} \d u \Bigg |<  
x^{-c} \int\limits_{\r} |\mm(c+\i u)| \d u =O(x^{-c})
\eeqq
which proves (\ref{kx_asympt_x_to_0}). The proof of (\ref{kx_asympt_x_to_infty}) is identical, except that we have to shift the contour of integration 
in the opposite direction.
\end{proof}

It turns out that for almost all parameters $\alpha$ the asymptotic series (\ref{kx_asympt_x_to_0}) and (\ref{kx_asympt_x_to_infty}) converge to $p(x)$ for all 
$x>0$. In order to state this result, we need to define the following set of real numbers.
\begin{definition}\label{def_set_L}
 Let ${\mathcal L}$ be the set of real irrational numbers $x$, for which there exists a constant $b>1$ such that the inequality
\beq\label{eqn_def_set_L}
 \bigg| x -\frac{p}{q} \bigg| < \frac{1}{b^{q}}
\eeq 
is satisfied for infinitely many integers $p$ and $q$.
\end{definition}

This set was introduced in \cite{Kuznetsov2010} in the connection with the distribution of the supremum of a stable process and it was later studied in \cite{KuzHub2010}. 
It was proved in \cite{KuzHub2010} that $x\in {\mathcal L}$ if and only if 
\beq\label{set_L_lim_condition}
 \lim\limits_{q\to +\infty} \frac{\ln \langle qx \rangle }{q} =0.
\eeq 
where $\langle x \rangle=\min\{ |x-i| \; : \; i\in {\mathbb Z}\}$. There also exists a characterization of elements of ${\mathcal L}$ in terms of their continued 
fraction expansion (see Proposition 1 in \cite{KuzHub2010}). 
It is known that ${\mathcal L}$ is a proper subset of Liouville numbers, that it is dense in $\r$ and that it has Lebesgue measure zero. 
The Hausdorff dimension of ${\mathcal L}$ is also zero. This set is closed under addition/multiplication by rational numbers. 
It is also known that $x\in {\mathcal L}$ if and only if $x^{-1} \in {\mathcal L}$.  
See \cite{KuzHub2010} for proofs of these results and for some further references.

The following Theorem is our main result in this section. 
\begin{theorem}\label{thm_density_expansion} Assume that $\alpha \notin {\mathcal L} \cup {\mathbb Q}$. Then for all $x>0$
 \beq\label{kx_series}
 p(x) = 
 \begin{cases}
  &\sum\limits_{n\ge 0} a_n x^n + \sum\limits_{m\ge 0} \sum\limits_{n\ge 0} b_{m,n} x^{(m+1-\beta+\gamma)\delta+n}, \; {\textnormal{ if }} \gamma+\hat \gamma<1, \\
  &\sum\limits_{m\ge 0} \sum\limits_{n\ge 0} c_{m,n} x^{-(m+\hat \beta)\delta-n-1},  \qquad\qquad\;\;\;\;\; {\textnormal{ if }} \gamma+\hat \gamma>1.
 \end{cases}
 \eeq
\end{theorem}

First let us establish the following technical result, which gives us a formula for $M(s)$ similar to the reflection formula for the Gamma function.
\begin{lemma}\label{lemma_mms_functional}
Define
\beq\label{def_ck}
c_k=1-(1-\beta+\gamma)\delta-\delta/2-k.
\eeq
Then for all $u \in \c$
 \beq\label{eqn_mms_functional}
 \mm(c_k+\i u)=(-1)^k\mm(c_k+k+\i u) \times \prod\limits_{j=0}^{k-1}
\frac{\cos(\pi \alpha(j-\i u + \gamma \delta))}{\cos(\pi \alpha(j  - \i u))} 
\times \frac{F(-\i u)}{F(-\i u + k)},
\eeq
where we have defined
\beq\label{def_Fw}
F(w)=\Gamma((1-\beta+\gamma)\delta+\delta/2+w)\frac{G(\frac{3}2\delta+w;\delta) }{G((\gamma+\frac{3}2)\delta+w;\delta)} \times
   \frac{G((\eta-\hat \gamma+\frac12)\delta+w;\delta)}
	 {G((\eta+\frac12)\delta+w;\delta)}.
\eeq
\end{lemma}
\begin{proof}
 Let us set $s=c_k+\i u$. Then iterating the functional identity (\ref{eqn_M_s_plus_1}) $k$ times we obtain
 \beq\label{eqn_mms1}
\mm(s)&=&\mm(s+k)\frac{\Gamma(s)}{\Gamma(s+k)} \prod\limits_{j=0}^{k-1}
 \frac{\Gamma(1-\beta+\gamma+\alpha s +\alpha j)}{\Gamma( 1-\beta+\alpha s +\alpha j)}
 \frac{\Gamma(\hat\beta+\hat\gamma-\alpha s -\alpha j)}{\Gamma( \hat\beta-\alpha s -\alpha j)}.
 \eeq
We use identity 
\beqq
 \frac{\Gamma(s)}{\Gamma(s+k)}=(-1)^k \frac{\Gamma(1-k-s)}{\Gamma(1-s)}
 \eeqq
and the reflection formula for the gamma function and rewrite (\ref{eqn_mms1}) as
\beq\label{eqn_mms2}
\mm(s)&=&(-1)^k\mm(s+k)\frac{\Gamma((1-\beta+\gamma)\delta+\delta/2-\i u)}{\Gamma(k+(1-\beta+\gamma)\delta+\delta/2-\i u)} \\
 &\times&
  \prod\limits_{j=0}^{k-1}
 \frac{\Gamma(-\alpha(k-1-j)-\frac12+\i \alpha u)}{\Gamma(-\alpha(k-1-j)-\gamma-\frac12+\i \alpha u)}
 \frac{\Gamma(\alpha(k-1-j)+\eta+\frac12-\i \alpha u)}{\Gamma(\alpha(k-1-j)+\eta-\hat \gamma+\frac12-\i \alpha u)}.
\eeq
Next we change the index $j\mapsto k-1-j$ and use reflection formula for the gamma function to obtain
\beq\label{eqn_mms3}
  \prod\limits_{j=0}^{k-1}
 \frac{\Gamma(-\alpha(k-1-j)-\frac12+\i \alpha u)}{\Gamma(-\alpha(k-1-j)-\gamma-\frac12+\i \alpha u)}= \prod\limits_{j=0}^{k-1}
\frac{\cos(\pi \alpha(j-\i u + \gamma \delta))}{\cos(\pi \alpha(j  - \i u))} 
\frac{\Gamma(\alpha j+\gamma+\frac32-\i \alpha u)}{\Gamma(\alpha j+\frac32-\i \alpha u)}.
\eeq
The proof of (\ref{eqn_mms_functional}) follows from (\ref{eqn_mms2}), (\ref{eqn_mms3}) and the following identity
\beqq
\prod\limits_{j=0}^{k-1} \Gamma(\alpha j + z)=\frac{G(\delta z+ k;\delta)}{G(\delta z;\delta)},
\eeqq
which is obtained by iterating formula (\ref{funct_rel_G_1}) $k$ times.
\end{proof}

{\it Proof of Theorem \ref{thm_density_expansion}:}
We will use a similar technique as in the proof of Theorem 2 in \cite{KuzHub2010}. 
Let us define $B=1-\gamma-\hat \gamma$ and assume that $B>0$. We start with (\ref{p_shift_contour}) and set $c=c_k$, where $c_k$ is defined by (\ref{def_ck}). 
Note that $\mm(s)$ does not have singularities on the vertical line $c_k+\i \r$; if this was not the case then $c_k$ would coincide with one of the poles $z^-_{m,n}$, which would imply that $\alpha$ is rational.

We define 
\beq\label{def_I1_I2}
I_1(x,k)=x^{-c_k}\re\left[\;\int\limits_{0}^{k} \mm(c_k+\i u) x^{-\i u} \d u \right], \;\;\; 
I_2(x,k)=x^{-c_k}\re\left[\;\int\limits_{k}^{\infty} \mm(c_k+\i u) x^{-\i u} \d u \right].
\eeq
It is clear that the integral in the right-hand side of (\ref{p_shift_contour}) is equal to $2(I_1(x,k)+I_2(x,k))$. Our goal is to prove 
that $I_j(x,k) \to 0$ as $k\to +\infty$ for all $x>0$. 

First, let us deal with $I_2(x,k)$. Using Theorem \ref{thm_main}, formula (\ref{eqn_Ms_asymptotics}) and Stirling's asymptotic formula (\ref{stirling}) 
  we find that there exists a constant $C_1>0$ such that
for all $s$ in the domain 
\beqq
 {\mathcal D}=\{s\in \c \; : \; |s|>2 {\textnormal{ and }}  \frac{\pi}8<\arg(s)<\frac{7\pi}8 \}
\eeqq
we have an upper bound
\beqq
|\mm(s)|< | s|^{C_1} \exp\left( \re\left[B s \ln(s) + s ((1+\ln(\delta))(1-B)+\pi \i \hat \gamma-1)\right]\right).  
\eeqq
From now on we will denote $u=\im(s)$. Computing the real part in the above expression we obtain that for all $s\in {\mathcal D}$
\beq\label{thm4_proof1}
|\mm(s)|&<& |s|^{C_1} \exp\left( B \re(s) \ln|s| -(\pi \hat \gamma+B\arg(s))u + \re(s) ((1+\ln(\delta))(1-B)-1) \right)
\eeq

Next we check that for all $k$ large enough, conditions $u>k$ and $s=c_k+\i u$ imply $s \in {\mathcal D}$, 
so the integrand in formula  (\ref{def_I1_I2}) defining $I_2(x,k)$ can be bounded from above as given in  (\ref{thm4_proof1}). 
Let us restrict $s$ to the line of integration $L_2=\{s=c_k+\i u \; :\; u\ge k\}$ and simplify this upper bound. Basically, we want to isolate a term of the 
form $\exp(- B k\ln(k))$ and show that everything else does not grow faster than an exponential function of $k$. First of all, 
when $s\in L_2$ we have $\re(s)=c_k<1-k$ and $|s|>|\re(s)|=|c_k|>k-1$, therefore
\beqq
\exp\left( B \re(s) \ln|s| \right)<\exp(- B (k-1)\ln(k-1)).
\eeqq
 Next, for $s\in L_2$ and $k$ sufficiently large it is true that $-2k<\re(s)<1-k$, which implies that there exists a constant $C_2>0$ such that 
\beqq
\exp\left( \re(s) ((1+\ln(\delta))(1-B)-1) \right)< C_2^k. 
\eeqq
Finally, for $s\in L_2$ we have $\arg(s)>0$, which together with the assumption $B>0$ shows that for all $u>0$
\beqq 
\exp\left(-(\pi \hat \gamma+B\arg(s))u\right)< \exp(-\pi \hat \gamma u).
\eeqq
Combining the above three estimates with (\ref{thm4_proof1}) and (\ref{def_I1_I2}) we see that
\beq
|I_2(x,k)| &\le & x^{-c_k} \int\limits_{k}^{\infty} |\mm(c_k+\i u)| \d u \\ \nonumber
&<& x^{-c_k}  C_2^k e^{-B (k-1)\ln(k-1)}
 \int\limits_{k}^{\infty} (c_k^2+u^2)^{\frac{C_1}2} e^{-\pi \hat \gamma u} \d u\\ \nonumber
&<& 
|c_k|^{C_1+1}  x^{-c_k}  C_2^k e^{-B (k-1)\ln(k-1)}.
\eeq
The right-hand side of the above inequality converges to $0$ as $k\to +\infty$, therefore $I_2(x,k)\to 0$ as $k\to +\infty$.

Now we will deal with $I_1(x,k)$. 
Our first goal is to find an upper bound for the product of trigonometric functions in 
(\ref{eqn_mms_functional}). We will follow the proof of Theorem 2 in \cite{KuzHub2010}: we use the trigonometric identities
\beqq
  \cos(x+\i y)&=&\cos(x)\cosh(y) - \i \sin(x) \sinh(y), \\
  | \cos(x+\i y)|^2 &=& \cosh(y)^2 - \sin(x)^2,
\eeqq
which imply that $|\cos(x) | \cosh(y) \le |\cos(x+\i y)| \le \cosh(y)$, therefore
\beqq
\bigg |\frac{\cos(a+\i y)}{\cos(b+\i y)} \bigg | \le \frac{1}{|\cos(b)|}.
\eeqq
Applying the above estimate and Lemma 1 from  \cite{KuzHub2010} we conclude that for $\alpha \notin  {\mathcal L} \cup {\mathbb Q}$ and for 
 $k$ large enough
\beqq
\Bigg | \prod\limits_{j=0}^{k-1}
\frac{\cos(\pi \alpha(j-\i u + \gamma \delta))}{\cos(\pi \alpha(j  - \i u))} \Bigg| \le 
\prod\limits_{j=0}^{k-1} |\sec(\pi \alpha j)|=2^{k+o(k)}<3^k. 
\eeqq
Using (\ref{def_I1_I2}), (\ref{eqn_mms1}) and the above inequality we conclude that for all $k$ large enough
\beq\label{thm4_proof2}
|I_1(x,k)|< x^{-c_k} 3^k \int\limits_{0}^{k} |\mm(c_k+k+\i u)| \times \bigg | \frac{F(-\i u)}{F(-\i u + k )} \bigg | \d u. 
\eeq

Now our goal is to prove that $F(\i u)/F(-\i u+k)$ converges to zero faster than any exponential function of $k$ as $k\to +\infty$. 
We use (\ref{def_Fw}), Stirling's formula (\ref{stirling}) and the asymptotic expansion (\ref{eqn_G2_asymptotics}) 
to conclude that when $w\to \infty$ in the domain $|\arg(w)|<3\pi/4$ we have
\beqq
\ln(F(w))=Bw \ln(w) + O(w).
\eeqq
This asymptotic result implies that there exists a constant $C_3>0$ such that for $k$ large enough and for all $v \in [0,1]$ 
\beqq
 \bigg | \frac{F(-\i k v)}{F(-\i k v + k )} \bigg |&<&
C_3^k \exp\left(B  \re\left[  (-\i k v) \ln(-\i k v) -  k(1-\i v) \ln(k(1-\i  v)) \right] \right)
\\
&=&
C_3^k \exp\left(B \left[ - k \ln (k\sqrt{1+v^2})+ k v \left(\arctan(v)-\frac{\pi}2 \right) \right]  \right)
\\ &<&
C_3^k \exp\left(-B k \ln(k)\right)
\eeqq
where in the last step we have used the fact that $\arctan(v)<\pi/2$ and $B>0$. Thus we see that for all $k$ large enough we have
\beqq 
\max\limits_{0\le u \le k } \bigg | \frac{F(-\i u)}{F(-\i u + k )} \bigg |=\max\limits_{0\le v \le 1} \bigg | \frac{F(-\i k v)}{F(-\i k  v + k )} \bigg |
< C_3^k \exp\left(-B k \ln(k)\right)
\eeqq
Combining the above estimate with (\ref{thm4_proof2}) we obtain that for all $k$ large enough
\beqq
|I_1(x,k)| < x^{-c_k} 3^k C_3^k e^{-B k \ln(k)} 
\int\limits_{0}^{\infty} |\mm(1-(1-\beta+\gamma)\delta-\delta/2+\i u)|  \d u 
\eeqq
and we see that $I_1(x,k)\to 0$ as $k\to +\infty$. Thus when $\gamma+\hat \gamma<1$ the first series in (\ref{kx_series}) converges to $p(x)$ for all $x>0$. 

When $\gamma+\hat \gamma>1$ the proof is very similar, except that one has to shift the contour of integration in the opposite direction. 
\qed

\section{Applications}\label{sec_applications}
 
In this section we  will present several applications of the above results on exponential functionals of the hypergeometric L\'evy processes. In particular, we 
will study various functionals related to strictly stable L\'evy processes. Our main tool will be the Lamperti transformation, 
which links a positive self-similar Markov process (pssMp) to an associated L\'evy process. By studying this associated L\'evy process and its exponential functional
we can obtain many interesting results about the original self-similar Markov process.

We assume that the reader is familiar with the Lamperti transformation \cite{LA} and Lamperti-stable processes (see \cite{CC, CKP, CPP, CPP1}). 

\subsection{Extrema of stable processes}\label{subsec_extrema_stable}
 
As our first example, we will present a new proof of some the results on the distribution of extrema of  general stable processes, which were obtained
in \cite{Kuznetsov2010} and \cite{KuzHub2010}. 
Let us define $S_1=\sup\{Y_t : \; 0\le t \le 1\}$, where $Y$ is a stable process with parameters $(\alpha,\rho)$, whose characteristic exponent is given by (\ref{def_Psi_stable}).
The infinite series representation for the density of $S_1$ (valid for almost all values of $\alpha$) was derived recently in \cite{KuzHub2010}. 
This formula was obtained using an explicit expression for the Mellin transform of $S_1$, which was found in \cite{Kuznetsov2010} via the Wiener-Hopf factorization of stable processes. Our goal in this section is to present a more direct route which leads to the density of $S_1$.

\begin{proposition}\label{prop_id_S1_I_alpha_xi}
 We have the following identity in distribution 
\beq\label{id_S1_I_alpha_xi}
S_1^{-\alpha}\stackrel{d}{=}I(\alpha,\xi),
\eeq
 where 
$\xi$ is a hypergeometric L\'evy process with parameters $(\alpha\rho, \alpha\rho, \alpha\rho, \alpha(1-\rho))$. 
\end{proposition}
\begin{proof}   Some of the  arguments in our proof have been used before in \cite{CKP}, nevertheless, we provide them here for completeness.  
Let us consider the dual process $\hat Y=-Y$. Note that the process $\hat Y$ is itself a stable 
process with parameters $(\alpha,1-\rho)$. Let $\hat Z$ be the process killed at the first exit from 
the positive half-line
\beqq
\hat Z_t=
\begin{cases}
& \hat Y_t, \;\;\; {\textnormal{ if }} t<\hat T, \\
& \Delta, \;\;\; {\textnormal{ if }} t\ge \hat T, 
\end{cases}
\eeqq 
where $\hat T=\inf\{t:\hat Y_t\le 0\}$ and $\Delta$ is the cemetery state. Let $\hat \xi$ be the process associated with $\hat Z$ by the Lamperti transformation.
 We use Theorem \ref{theorem_hyper=Lamperti_stable} and find that $\hat \xi$ is a hypergeometric process with parameters $(1-\alpha \rho,\alpha(1-\rho),1-\alpha \rho,\alpha \rho)$. 

Next, from the Lamperti transformation we find that $\hat T\stackrel{d}{=}I(\alpha,-\hat\xi)$. This statement is equivalent 
to $\hat T\stackrel{d}{=}I(\alpha,\xi)$, where $\xi=-\hat \xi$ is a hypergeometric process with parameters $(\alpha \rho,\alpha\rho,\alpha \rho,\alpha (1-\rho))$
(see Proposition \ref{prop1}(v)). The rest of the proof follows from the chain of identities 
\beqq
\p_1(\hat T>t)&=&\p_1(\inf_{0\leq s\leq t}\hat Y_s>0)=\p_0(\inf_{0\leq s\leq t}\hat Y_s>-1)\\
&=&\p_0(\sup_{0\leq s\leq t} Y_s<1)=\p_0(t^{1/\alpha}\sup_{0\leq s\leq 1} Y_s<1)=\p_0(S_1^{-\alpha}>t).
\eeqq
 \end{proof}

Proposition \ref{prop_id_S1_I_alpha_xi}  implies that the density of $S_1$ can be represented in terms of the density of $I(\alpha,\xi)$ (denoted by $p(x)$) as follows
\beqq
\frac{\d}{\d x} \p(S_1\le x)=\alpha x^{-1-\alpha}p(x^{-\alpha}), \qquad x>0.
\eeqq
Using this expression and Theorem \ref{thm_density_asymptotics} \{Theorem \ref{thm_density_expansion}\} 
 we recover the asymptotic expansions that appears in Theorem 9 in \cite{Kuznetsov2010}  \{the convergent series representations given in Theorem 2 in \cite{KuzHub2010}\}.\\

{\bf Remark 2:}
Note that the convergent and asymptotic series given in Theorems  \ref{thm_density_expansion} and \ref{thm_density_asymptotics} 
have identical form. Combining these two results we obtain the following picture: depending on whether $\gamma+\hat \gamma<1$ or $\gamma+\hat \gamma>1$ one of the 
series in (\ref{kx_series}) converges for all $x \in (0,\infty)$, and this convergent series is also asymptotic at one of the boundaries 
of this interval (at $0$ or $+\infty$). At the same time the second 
series  in (\ref{kx_series}) provides an asymptotic series at the other boundary. Therefore in what follows we will only present 
the statements about the convergent series representation, with the understanding that the asymptotic expansions are implicitly imbedded into these results. 

\subsection{Entrance laws}\label{subsec_entrance_laws}

In this section we will obtain some explicit results on the entrance law of a stable process conditioned to stay positive and the entrance 
law of the excursion measure of the stable process reflected at its infimum. 
We will use the following result from \cite{BeY}: if a non-arithmetic L\'evy process $X$ satisfies
$\e[|X_1|]<\infty$ and $\e[X_1]>0$, then as $x\to 0+$ its corresponding
pssMp $(Z,\mathbf{Q}_{x})$ in the Lamperti representation converges weakly
(in the sense of finite-dimensional distributions) towards a nondegenerate probability law $(Z,\mathbf{Q}_0)$.
Under these conditions, the entrance law under $\mathbf{Q}_0$ is described
as follows: for every $t>0$ and every measurable function $f:\r_+\to
\r_+$,
\begin{equation}\label{entlaw}
\mathbf{Q}_0\big[f(Z_t)\big]=\frac{1}{\alpha \mathbb{E}[X_1]}\e\left[\frac{1}{I}
f\left(\frac{t}{I}\right)\right]\,,
\end{equation}
where $I=I(\alpha,X)=\int_{0}^{\infty} \exp\{-\alpha X_{s}\}\,{\rm d} s.$ Necessary and sufficient conditions for the weak  convergence of $(Z, \mathbf{Q}_x)$ on  the Skorokhod space were given in
\cite{CCh}. 

As before, we consider a stable process $(Y,\p_x)$ and denote by $(Y, \p_x^{\uparrow})$ the stable process conditioned to stay positive 
(see \cite{Ch,CD} for a proper definition). It is known that $(Y, \p_x^{\uparrow})$ is a pssMp with index $\alpha$. 
According to Corollary 2 in \cite{CC} (and Theorem \ref{theorem_hyper=Lamperti_stable} above) its associated L\'evy process is given by $\xi^\uparrow$, which is 
a hypergeometric L\'evy process with parameters $(1,\alpha\rho, 1,\alpha(1-\rho))$.
From Proposition \ref{prop1}(ii) we see that $\e[\xi^{\uparrow}_1]>0$, and the fact that hypergeometric processes have L\'evy measures with exponentially decaying tails implies that $\e[|\xi^{\uparrow}_1|]<\infty$, thus we have the weak convergence of  $(Y, \p_x^{\uparrow})$ as $x\to 0^+$. 
We denote the limiting law by $\p^{\uparrow}$. Note also that in this particular case, the weak convergence of $(Y, \p_x^{\uparrow})$ has been proved in a direct way in
\cite{Ch}. 

 Everywhere in this section the coefficients $ \{b_{m,n}\}_{m,n\ge 0}$ and $\{c_{m,n}\}_{m,n\ge 0}$ are defined as in Definition \ref{def_an_bmn_cmn} with  parameters 
 $\beta=\hat{\beta}=1$, $\gamma=\alpha\rho$ and $\hat{\gamma}=\alpha(1-\rho)$. 
\begin{proposition}\label{prop_entrance_laws}
Let $p^\uparrow_t$ be the density of  the entrance law of 
$(Y,\p^{\uparrow})$. If $\alpha\notin\mathcal{L}\cup \mathbb{Q}$ then for $x>0$ we have 
\beq\label{eqn_p1_up}
 p_1^\uparrow(x) = 
 \begin{cases}
  &\frac{x^{-1-\rho}}{\alpha \Gamma(\alpha \rho) \Gamma(1+\alpha(1-\rho))} \sum\limits_{m\ge 0} \sum\limits_{n\ge 0} b_{m,n} x^{-m/\alpha-n},\;\;\; {\textnormal{ if }} \alpha<1, \\
  &\frac{x^{1/\alpha}}{\alpha \Gamma(\alpha \rho) \Gamma(1+\alpha(1-\rho))}\sum\limits_{m\ge 0} \sum\limits_{n\ge 0} c_{m,n} x^{m/\alpha+n},  \;\;\; \;\;{\textnormal{ if }} \alpha>1.
 \end{cases}
 \eeq
\end{proposition}
\begin{proof}
Let $p(x)$ be the density of $I(\alpha,\xi^{\uparrow})$. From (\ref{entlaw}) we find that the density of the entrance law of  $(Y,\p^{\uparrow})$
is given by
\beqq
p_1^\uparrow(x)=\frac{1}{\alpha \e[\xi_1^{\uparrow}]} x^{-1}p(x^{-1}).
\eeqq
Using Proposition \ref{prop1}(ii) we obtain $\e[\xi^{\uparrow}_1]=\Gamma(\alpha \rho)\Gamma(1+\alpha(1-\rho))$. In order to finish the proof we only need
to apply results of Theorem \ref{kx_series}. 
\end{proof}

\begin{proposition}\label{prop_excursion_entrance_law}
  Let $q_t$ be the density of  the entrance law of the excursion measure of the reflected process $(Y-\underline{Y}, \p)$, where 
$\underline{Y}_t=\inf_{0\leq s\leq t}Y_s$. If $\alpha\notin\mathcal{L}\cup \mathbb{Q}$ then for $x>0$ we have 
\beq\label{eqn_q1}
 q_1(x) = 
 \begin{cases}
  &\frac{x^{-1-\rho-\alpha(1-\rho)}}{\alpha \Gamma(\alpha \rho) \Gamma(1+\alpha(1-\rho))} \sum\limits_{m\ge 0} \sum\limits_{n\ge 0} b_{m,n} x^{-m/\alpha-n},\;\;\; {\textnormal{ if }} \alpha<1, \\
  &\frac{x^{1/\alpha-\alpha(1-\rho)}}{\alpha \Gamma(\alpha \rho) \Gamma(1+\alpha(1-\rho))}\sum\limits_{m\ge 0} \sum\limits_{n\ge 0} c_{m,n} x^{m/\alpha+n},  \;\;\; \;\;{\textnormal{ if }} \alpha>1.
 \end{cases}
 \eeq
\end{proposition}
\begin{proof}
The proof follows from Proposition \ref{prop_entrance_laws} and Theorem 3 in \cite{Ch}, which gives us the identity $q_1=x^{-\alpha(1-\rho)}p_1^{\uparrow}(x)$.
\end{proof}

\subsection{The distribution of the lifetime of a stable process conditioned to hit zero continuously}\label{subsec_continuously_zero}

Our third application deals with stable
processes conditioned to hit zero continuously. Processes in this class are defined as Doob's $h$-transform with respect to the function $h(x)=
x^{\alpha(1-\rho)-1}$, which is excessive  for the killed process
$(Y_t\mathbf{1}_{\{t<T\}},\p_x)$. Its law $\p^\downarrow_x$,
which is defined on each $\sigma$-field $\mathcal{F}_t$ by
\begin{equation}\label{downarrow}
\left. \frac{{\rm d} \mathbb{P}^{\downarrow}_x}{{\rm d} \mathbb{P}_x}\right|_{\mathcal{F}_t}=\frac{Y^{\alpha(1-\rho)-1}_t}{x^{\alpha(1-\rho)-1}}\mathbf{1}_{\{t<T\}},
\end{equation} 
 is that of a pssMp that hits zero in a continuous way. According to Corollary 3 in \cite{CC} (and Theorem \ref{theorem_hyper=Lamperti_stable} above), 
its associated L\'evy process in the Lamperti
representation   is the hypergeometric L\'evy process $\xi^\downarrow$ with parameters $(0,\alpha\rho,0,\alpha(1-\rho))$. 
In particular from Proposition \ref{prop1}(ii) we find that the process $\xi^{\downarrow}$ drifts to $-\infty$. Let $T^\downarrow$ be the life-time of  the stable
process conditioned to hit zero continuously. From the Lamperti transformation it follows that $T^\downarrow\stackrel{d}{=}I(\alpha, -\xi^\downarrow)$, therefore
$T^\downarrow\stackrel{d}{=}I(\alpha, \xi)$, where $\xi$ is a hypergeometric process with parameters $(1,\alpha(1-\rho),1,\alpha\rho)$.
Let us define the coefficients $ \{b_{m,n}\}_{m,n\ge 0}$ and $\{c_{m,n}\}_{m,n\ge 0}$  as in Definition \ref{def_an_bmn_cmn} with  $\beta=\hat{\beta}=1$, $\gamma=\alpha(1-\rho)$ and $\hat{\gamma}=\alpha\rho$. Then according to Theorem \ref{kx_series}, if $\alpha\notin \mathcal{L}\cup \mathbb{Q}$ and $x>0$ we have  
 \beqq
 \frac{\d}{\d t} \p_1^\downarrow ( T^\downarrow \le t) = 
 \begin{cases}
  & t^{1-\rho} \sum\limits_{m\ge 0} \sum\limits_{n\ge 0} b_{m,n} t^{m/\alpha+n},\qquad\;\;\; {\textnormal{ if }} \alpha<1, \\
  & t^{-1-1/\alpha}\sum\limits_{m\ge 0} \sum\limits_{n\ge 0} c_{m,n} t^{-m/\alpha-n},  \;\;\; {\textnormal{ if }} \alpha>1.
 \end{cases}
 \eeqq

\subsection{Distribution of some functionals related to the radial part of a symmetric stable process}\label{subsec_radial_part}
 
Our last application deals with the radial part of symmetric stable processes in $\mathbb{R}^d$. Let ${\bf Y}=({\bf Y}_t, t\geq 0)$ 
be a symmetric stable  L\'evy process of index $\alpha\in (0,2)$ in $\mathbb{R}^d$ ($d\ge 1$), defined by
\[
\e_0\big[e^{i \langle \lambda, {\bf Y}_t \rangle }\big]=e^{-t\|\lambda\|^{\alpha}},
\]
for all $t\ge 0$ and $\lambda \in \mathbb{R}^d$. Here $\p_y$ denotes the law of the process ${\bf Y}$ started  from $y\in \mathbb{R}^d$,  $\|\cdot\|$ is the norm in  $\mathbb{R}^d$ and $\langle\cdot,\cdot\rangle$ is the Euclidean inner product. 

According to Caballero et al. \cite{CPP1}, when $\alpha < d$  the radial process $R_t=\frac12\|{\bf Y}_t\|$ is a transient 
positive self-similar Markov process with index $\alpha$ 
and infinite lifetime. From Theorem 7 in \cite{CPP1} we find that the Laplace exponent of its associated L\'evy process $\xi$  is given by 
\begin{equation}\label{WHSS}
\psi_{\xi}(z)=-\frac{\Gamma((\alpha-z)/2)}{\Gamma(-z/2)}\frac{\Gamma((z +d)/2)}{\Gamma((z +d-\alpha)/2)}.
\end{equation}
This shows that  the process $X:=2\xi$ is a hypergeometric L\'evy process with parameters $(1,\alpha/2,(d-\alpha)/2,\alpha/2)$. From Proposition \ref{prop1}(ii) 
we see that $\xi$ drifts to $+\infty$ and that 
\beq
\e[\xi_1]=\frac{1}{2}\frac{\Gamma\left(\frac{\alpha}2\right) \Gamma\left(\frac{d}2\right)}{\Gamma\left(\frac{d-\alpha}2\right)}.
\eeq 

Let $\mathbf{P}_x$ denotes the law of the process $R$ starting from $x>0$. The process $\xi$ satisfies $\e[\xi_1]>0$ and $\e[|\xi_1|<\infty$, and according 
to \cite{BeY, CCh} we have the weak convergence of $(R, \mathbf{P}_x)$ towards $(R,\mathbf{P}_0)$ as $x\to 0^+$.

\begin{proposition}\label{prop_radial_1} Let $\tilde p_t$ be the density of  the entrance law of 
$(R,\mathbf{P}_0)$. If $\alpha<d$ then for $x>0$ we have 
\beq\label{eqn_radial_1}
 \tilde p_1(x) = 
 \begin{cases}
  &\frac{2x^{-1}}{\pi \alpha \Gamma\left(\frac{d}{2}\right)}  \sum\limits_{n\ge 1} \sin\left(\frac{\pi \alpha n}2\right)
 \Gamma\left(1+\frac{\alpha n}2\right) \Gamma\left(\frac{d+\alpha n}2\right) \frac{(-1)^{n+1}}{n!}   x^{-n},  \;\;\; {\textnormal{ if }} \alpha<1, \\
  &\frac{4x^{-1+d/\alpha}}{\alpha^2  \Gamma\left(\frac{d}{2}\right)}  
\sum\limits_{m\ge 0} \frac{\Gamma\left(\frac{d+2 m}{\alpha}\right)}{ \Gamma\left(\frac{d+2 m}2\right)} \frac{(-1)^{m}}{m!}   x^{2m/\alpha}
\qquad\qquad\qquad\qquad\;\;\;\;\;\; {\textnormal{ if }} \alpha>1.
 \end{cases}
 \eeq
\end{proposition}
\begin{proof}
We follow the same approach as in the proof of  Proposition \ref{prop_entrance_laws}. Using (\ref{entlaw}) we conclude that
$\tilde p_1(x)=(2/\alpha \e[\xi_1]) x^{-1}p(x^{-1})$, where $p(x)$ is the density of $I(\alpha,\xi)=I(\alpha/2,X)$. In the case 
$\alpha\notin \mathcal{L}\cup \mathbb{Q}$ the infinite 
series representation for $p(x)$ and the expression for the coefficients follow from Theorem \ref{kx_series} and Definition 
\ref{def_an_bmn_cmn} with $\beta=1$, $\gamma=\alpha/2$, $\hat \beta = (d-\alpha)/2$ and $\hat \gamma=\alpha/2$. Note that both infinite series
in (\ref{eqn_radial_1}) converge for all $\alpha$, thus we can remove the assumption $\alpha\notin \mathcal{L}\cup \mathbb{Q}$.

There is also a simpler way to derive an explicit series representation for $p(x)$. Theorem \ref{thm_main} and Definition \ref{def_Ms} tell us that the Mellin transform of $I(\alpha/2,X)$ is given by
\beq\label{mellin_radial}
\e[I(\alpha/2,X)^{s-1}]=\Gamma(s) M\left(s;\frac{\alpha}2,1,\frac{\alpha}2,\frac{d-\alpha}2,\frac{\alpha}2\right)=
C \Gamma(s) \frac{G \left(s; \frac{2}{\alpha}\right)}{G \left(1+s; \frac{2}{\alpha}\right) }
 \frac{G \left(\frac{d}{\alpha}+1-s; \frac{2}{\alpha}\right) }{G \left(\frac{d}{\alpha}-s; \frac{2}{\alpha}\right) },
\eeq
where $C$ is a constant such that the right-hand side of the above identity is equal to one for $s=1$. Using the functional identity (\ref{funct_rel_G_1}) for the double
gamma function we simplify the right-hand side in (\ref{mellin_radial}) and obtain
\beq\label{mellin_alpha/2}
\e[I(\alpha/2,X)^{s-1}]=\frac{\Gamma\left(\frac{\alpha}2\right)}{\Gamma\left(\frac{d-\alpha}2\right)}\Gamma(s)\frac{\Gamma\left(\frac{d-\alpha s}2\right)}{\Gamma\left(\frac{\alpha s}2\right)}.
\eeq
We see that $\e[I(\alpha/2,X)^{s-1}]$ has simple poles at points $\{-n\}_{n\ge 1}$ and $\{(d+2m)/\alpha\}_{m\ge 0}$. The residues at these points (which give us the 
coefficients in (\ref{eqn_radial_1})) can be easily found using 
the fact that the residue of $\Gamma(s)$ at $s=-n$ is equal to $(-1)^n/n!$.
\end{proof}

Next, we will study the last passage time of $({\bf Y},\p_0)$ from the sphere in $\mathbb{R}^d$ of radius $r$, i.e.
\beqq
U_r=\sup\{s\ge 0: \|{\bf Y}_s\|<r\}.
\eeqq
Due to the self-similarity property of $Y$ we find that random variables $U_r$ satisfy the scaling property $b^{\alpha}U_{r}\stackrel{d}{=} U_{br}$ (valid 
for any $b,r>0$), thus it is sufficient to consider the case of of $U_2$. 
\begin{proposition}\label{prop_radial_2} If $\alpha<d$ then
\beq\label{eqn_radial_2}
 \frac{\d}{\d t} \p(U_2 \le t) = 
 \begin{cases}
  &\frac{1}{\pi  \Gamma\left(\frac{d-\alpha}{2}\right)}  \sum\limits_{n\ge 0} \sin\left(\frac{\pi \alpha (1+ n)}2\right)
 \Gamma\left(\frac{\alpha (1+n)}2\right) \Gamma\left(\frac{d+\alpha n}2\right) \frac{(-1)^{n}}{n!}   t^{n},  \;\;\; {\textnormal{ if }} \alpha<1, \\
  &\frac{2t^{-d/\alpha}}{\alpha  \Gamma\left(\frac{d-\alpha}{2}\right)}  
\sum\limits_{m\ge 0} \frac{\Gamma\left(\frac{d+2 m}{\alpha}\right)}{ \Gamma\left(\frac{d-\alpha}2+m+1\right)} \frac{(-1)^{m}}{m!}   t^{-2m/\alpha},
\qquad\qquad\qquad\qquad {\textnormal{ if }} \alpha>1.
 \end{cases}
 \eeq
\end{proposition}
\begin{proof}
Let us define the last passage time of $R$ as $L_r=\sup\{s\ge 0: \|R_s\|<r\}$. It is clear that $U_2\stackrel{d}{=}L_1$. 
According to Proposition 1 in \cite{CP}, the random variable $L_1$ has the same law as ${\mathcal G}^{\alpha} I(\alpha, \xi)$, where $\xi$ is the L\'evy process associated to $R$ by the Lamperti transformation and ${\mathcal G}$ is independent of $\xi$. Lemma 1 in \cite{CP}  tells us that 
${\mathcal G}\stackrel{d}{=}e^{-\mathcal{U}\mathcal{Z}}$, where $\mathcal{U}$ and $\mathcal{Z}$ are independent random variables, such that $\mathcal{U}$ is uniformly distributed over $[0,1]$ and the law of $\mathcal{Z}$ is given by
\beqq
\p(\mathcal{Z}>u)=\frac{1}{\e[H_1]} \int\limits_u^{\infty}  x\,\nu(\d x), \qquad u\ge 0,
\eeqq
where $H=(H_t,t\ge 0)$ denotes the ascending ladder height of $\xi$ and $\nu$ its L\'evy measure. From the proof of Lemma 3 in \cite{CPP1} we know that $H$ has no linear drift and that 
 its L\'evy measure is given by (up to a multiplicative constant)
\beqq
\nu(\d x)=\frac{e^{2x}}{(e^{2x}-1)^{1+\alpha/2}}\d x.
\eeqq
Using the above expression, integration by parts, the integral identity
(\ref{beta_integal}) and reflection formula for the gamma function we obtain  
\beqq
\e[H_1]=\int\limits_0^{\infty} x \nu(\d x)=\frac{\pi}{\alpha \sin(\frac{\pi \alpha}2)}.
\eeqq

Let us find the Mellin transform of $L_1$. We use the  independence of $U$ and ${\mathcal Z}$ and obtain
\beqq
\e[{\mathcal G}^{s}]&=&\e\left[e^{-s {\mathcal U} {\mathcal Z}}\right]=\frac{1}{s}\e \left[ \frac{1-e^{-s {\mathcal Z}}}{{\mathcal Z}} \right]=
\frac1{s \e[H_1]} \int\limits_0^{\infty}\frac{(1-e^{-s z})e^{2z} }{(e^{2z}-1)^{1+\alpha/2}}\d z\\
&=&-\frac{1}{\alpha s \e[H_1]} \int\limits_0^{\infty} (1-e^{-s z}) \d \left[ (e^{2z}-1)^{-\alpha/2}\right] \\
&=&\frac{1}{\alpha \e[H_1]}\int\limits_0^{\infty} \frac{e^{-s z}}{(e^{2z}-1)^{\alpha/2}} \d z
=\frac{1}{\alpha \e[H_1]} \frac{\Gamma\left(\frac{\alpha}2(1+s)\right)\Gamma\left(1-\frac{\alpha}2\right)}{\Gamma\left(1+\frac{\alpha s}2\right)},
\eeqq
where we have again applied integration by parts and integral identity (\ref{beta_integal}). Combining the above two expressions with 
(\ref{mellin_alpha/2}) we have
\beqq
\e[L_1^{s-1}]=\e[{\mathcal G}^{s-1}]\times \e[I(\alpha,\xi)^{s-1}]=\frac{\Gamma(s)}{\Gamma\left(\frac{d-\alpha}2\right)}\frac{\Gamma\left(\frac{d-\alpha s}2\right)}{\Gamma\left(1-\frac{\alpha(1-s)}2\right)}.
\eeqq
The function in the right-hand side of the above equation has simple poles at points $\{-n\}_{n\ge 0}$ and $\{(d+2m)/\alpha\}_{m\ge 0}$. 
We express the density of $L_1$ as the inverse Mellin transform, compute the residues and 
use similar technique as in the proof of Theorem \ref{kx_series} to obtain the series representation (\ref{eqn_radial_2}).
\end{proof}

\newpage


\end{document}